\magnification=\magstep1
\vbadness=10000
\hbadness=10000
\tolerance=10000

\proclaim Quantum vertex algebras.  \hfill 7 March 1999. 

Richard E. Borcherds, 
\footnote{$^*$}{ Supported by a Royal Society
professorship. This paper was written at the Max-Planck institute in Bonn. }

D.P.M.M.S.,
16 Mill Lane, 
Cambridge, 
CB2 1SB,
England.

e-mail: reb@dpmms.cam.ac.uk

home page:  www.dpmms.cam.ac.uk/\hbox{\~{}}reb

\bigskip

\proclaim Contents. 

1. Introduction.

Notation. 

2. Twisted group rings. 

3. Construction of some categories. 

4. Examples of vertex algebras. 

5. Open problems.

\proclaim
1.~Introduction.

The purpose of this paper is to make the theory of vertex algebras
trivial.  We do this by setting up some categorical machinery so that
vertex algebras are just ``singular commutative rings'' in a certain
category.  This makes it easy to construct many examples of vertex
algebras, in particular by using an analogue of the construction of a
twisted group ring from a bicharacter of a group. We  also define
quantum vertex algebras as singular braided rings in the same
category and construct some examples of them. The constructions work
just as well for higher dimensional analogues of vertex algebras,
which have the same relation to higher dimensional quantum field
theories that vertex algebras have to one dimensional quantum field
theories.

One way of thinking about vertex algebras is to regard them as
commutative rings with some sort of singularities in their
multiplication.  In algebraic geometry there are two sorts of
morphisms: regular maps that are defined everywhere, and rational maps
that are not defined everywhere. It is useful to think of a
commutative ring $R$ as having a regular multiplication map from
$R\times R$ to $R$, while vertex algebras only have some sort of
rational or singular multiplication map from $R\times R$ to $R$ which
is not defined everywhere. One of the aims of this paper is to make
sense of this, by defining a category whose bilinear maps can be
thought of as some sort of maps with singularities.

The main idea for constructing examples of vertex algebras in this paper
is a generalization of the following well-known method for
constructing twisted group rings from bicharacters of groups. 
Suppose that $L$ is a discrete group (or monoid)
and $R$ is a  commutative ring. Recall that an $R$-valued bicharacter of $L$
is a map $r:L\times L\mapsto R$ such that
$$\eqalign{
r(1,a)&= r(a,1)=1\cr
r(ab,c)&=r(a,c)r(b,c)\cr
r(a,bc)&=r(a,b)r(a,c).\cr
}$$
If $r$ is any $R$-valued 
bicharacter of $L$ then we define a new associative multiplication 
$\circ$ on
the group ring $R[L]$ by putting
$a\circ b=abr(a,b)$. We call $R[L]$ with this new multiplication 
 the twisted group ring of $L$.
The point is that this rather trivial construction can be generalized 
from group rings to  bialgebras in additive symmetric tensor categories. 
We will construct vertex algebras by applying this construction
to ``singular bicharacters'' of bialgebras 
in a suitable additive symmetric tensor category. 

Section 2 describes how to generalize the twisted group ring construction
to bialgebras, and constructs several examples of singular bicharacters
that we will use later. Much of section 2 uses an extra structure
on the spaces underlying many common vertex algebras that is
often overlooked. It is well known that these spaces often 
have natural ring structures, but what is less well known is that
this can usually be extended to a cocommutative bialgebra structure. 
The comultiplication turns out to be very useful for keeping track 
of the behavior of vertex operators; this is not so important
for vertex algebras, but is very useful for quantum vertex algebras.
It also allows us to interpret
these spaces as the coordinate rings of gauge groups. 

Section 3 contains most of the hard work of this paper. We have to
construct a category in which the commutative rings are more or less
the same as vertex algebras. The motivation for the construction of
this category comes from classical and quantum field theory (though it
is not necessary to know any field theory to follow the
construction). The idea is to construct categories which capture all
the formal operations one can do with fields.  For examples, fields
can be added, multiplied, differentiated, multiplied by functions on
spacetime, and we can change variables and restrict fields.
All of these operations are trivial but there are so many of them 
that it takes some effort to write down all the compatibility 
conditions between them. The categories constructed in section
3 are really just a way of writing down all these compatibility
conditions explicitly. 
The main point of doing this is the definition at the end of section 
3, where we define $(A,H,S)$ vertex algebras to be the commutative rings
in these categories. Here $A$ is a suitable additive category (for example
the category of modules over a commutative ring), $H$ is a suitable
bialgebra in $A$ (and can be thought of as a sort of group ring
of the group of automorphisms of spacetime), and $S$ is something
that controls the sort of singularities we allow.  

One of the main differences between the $(A,H,S)$ vertex algebras
defined in section 3 and previous definitions is as follows. Vertex
algebras as usually defined consist of a space $V(1)$ with some extra
operations, whose elements can be thought of a fields depending one
one spacetime variable. On the other hand $(A,H,S)$ vertex algebras
include spaces $V(1,2,\ldots, n)$ which can be thought of as fields
depending on $n$ spacetime variables for all $n$. The lack of these
fields in several variables seems to be one reason why classical
vertex algebras are so hard to handle: it is necessary to 
reconstruct these fields, and there seems to be no canonical way to do
this. However if these fields are given in advance then a lot of these
technical problems just disappear.

Section 4 puts everything together to construct many examples of
vertex algebras.
The main theorem of this paper is theorem 4.2, which shows how to
construct a vertex algebra from a singular bicharacter of a
commutative and cocommutative bialgebra. As examples, we show that
the usual vertex algebra of an even lattice can be constructed like
this from the Hopf algebra of a multiplicative algebraic group, and
the vertex algebra of a (generalized) free quantum field theory can be
constructed in the same way from the Hopf algebra of an additive
algebraic group.  (This shows that the vertex algebra of a lattice is
in some sense very close to a free quantum field theory: they have the
same relation as multiplicative and additive algebraic groups.)

The vertex algebras we construct in this paper do not at first sight
look much like classical vertex algebras: they seem to be missing
all the structure such as vertex operators, formal power series, 
contour integration,
operator product expansions, and so on. We show that all this extra structure
can be reconstructed from the more elementary operations
we provide for vertex algebras. For example, the usual locality 
property of vertex operators follows from the fact that
we define vertex algebras as {\sl commutative} rings in some category. 

All the machinery in sections 2 and 3 has been set up so that it
generalizes trivially to quantum vertex algebras and higher
dimensional analogues of vertex algebras. For example, we define
quantum vertex algebras to be braided (rather than commutative) rings
in a certain category, and we can instantly construct many examples of
them from non-symmetric bicharacters of bialgebras.  By changing a
certain  bialgebra $H$ in the construction, we immediately get the
``vertex $G$ algebras'' of [B98], which have the same relation to
higher dimensional quantum field theories that vertex algebras have to
one dimensional quantum field theories.

Finally in section 5 we list some open problems and topics
for further research.

Some related papers are [F-R] and [E-K], which give alternative
definitions of quantum vertex algebras. These definitions are not
equivalent to the ones in this paper, but define concepts that are
closely related (at least in the case of 1 dimensional spacetime) in
the sense that the interesting examples for all definitions should
correspond.  Soibelman has introduced other foundations
for quantum vertex algebras, which seem to be related to 
this paper. There is also a preprint [B-D] which defines vertex
algebras as commutative rings or Lie algebras in suitable multilinear
categories. (Soibelman pointed out to me that multi categories seem
to have been first introduced by Lambek in [L].) 
It might be an interesting question to study the relationship of
this paper to [B-D].  One major difference is that the
paper [B-D] extends the
genus 0 Riemann surfaces that appear in vertex algebra theory to
higher genus Riemann surfaces, while in this paper we extend them
instead to higher dimensional groups.

I would like to thank S. Bloch, I. Grojnowski, 
J. M. E. Hyland,   and Y. Soibelman  for their help. 

\proclaim Notation.

\item{$A$} 
An additive symmetric tensor category.
\item{$C$} 
A symmetric tensor category, with tensor product $\cup$; usually $Fin$ 
or $Fin^{\not\equiv}$. 
\item{$\Delta$} 
The coproduct of a  bialgebra, or a propagator. 
\item{$D^{(i)}$}
An element of the formal group ring of the one dimensional additive
formal group.
\item{$\eta$} 
The counit of a bialgebra.
\item{$Fin$} 
The category of finite sets. 
\item{$Fin^{\not\equiv}$}
The category of finite sets with an inequivalence relation. 
\item{$Fun$}
A functor category. 
\item{$H$}
A cocommutative bialgebra in $A$.
\item{$I,J$}
Finite sets. 
\item{$L$}
An integral lattice.
\item{$M$}
A commutative ring in $A$, or a commutative cocommutative bialgebra. 
\item{$r$}
A bicharacter.
\item{$R$}
A commutative ring or an $R$-matrix.
\item{$R\hbox{-mod}$}
The category of $R$-modules.
\item{$S$}
A commutative ring in some category, especially $Fun(C,A, T^*(H))$.
\item{$S^*$}
A symmetric algebra. 
\item{$T$} 
A cocommutative Hopf algebra in $Fun(C^{op}, A)$. 
\item{$T_*$, $T^*$}
$T_*(M)(I)= \otimes_{i\in I} M$, $T^*(H)(I)= \otimes_{i\in I} H$. 
See definition 3.3.
\item{$U, V$}
Objects of $Fun(C, A, T^*(H), S)$. 

\proclaim 2.~Twisted group rings. 

We let $R$ be any commutative ring. Recall that a bialgebra
is an algebra with a compatible coalgebra structure,
with the coproduct and counit denoted by $\Delta$ and $\eta$,  and a Hopf
algebra is a bialgebra with an antipode. If $a$ is an element of
a coalgebra then we put $\Delta(a)=\sum a'\otimes a''$. 

Recall from the introduction that any bicharacter $r$ of a group $L$
can be used to define a twisted group ring. 
We now extend this idea from group rings $R[L]$ to cocommutative bialgebras. 

\proclaim Definition 2.1. Suppose that $M$ and $N$ are bialgebras over
$R$ and $S$ is a commutative $R$-algebra. Then we define a bimultiplicative
map 
from $M\otimes N$ to $S$
to be a linear map  $r:M\otimes N\mapsto S$ such that
$$\eqalign{
r(1\otimes a)&= \eta(a), \qquad r(a\times 1)=\eta(a)\cr
r(ab\otimes c)&=\sum r(a\otimes c')r(b\otimes c'')\cr
r(a\otimes bc)&=\sum r(a'\otimes b)r(a''\otimes c)\cr
}
$$
where $\Delta(a)=\sum a'\otimes a''$, $\Delta(c)=\sum c'\otimes c''$,
and $\eta$ is the counit of $M$ or $N$. We define an $S$-valued bicharacter
of $M$ to be a bimultiplicative map from $M\otimes M$ to $S$. 
We say the bicharacter $r$ is symmetric if $r(a\otimes b)=r(b\otimes
a)$ for all $a,b\in M$.

The $S$-valued bicharacters form a  monoid, which is commutative
if $M$ is co-commutative. The identity 
bicharacter is defined by $r(a\otimes b)=\eta(a)\otimes \eta(b)$, 
and the product $rs$ of two bicharacters $r$ and $s$ is given by
$$rs(a\otimes b)=\sum r(a'\otimes b')s(a''\otimes b'').$$
If $M$ is a Hopf algebra with antipode $s$ then any $S$-valued
bicharacter $r$ has an inverse $r^{-1}$ defined by
$$r^{-1}(a\otimes b)=r(s(a)\otimes b)$$
so the $S$-valued bicharacters form a group. 

\noindent {\bf Example 2.2.}
Suppose that $M=R[L]$ is the group ring of a group $L$, considered
as a bialgebra in the usual way (with $\Delta(a)=a\otimes a$ for $a\in L$).
Any $R^*$-valued bicharacter of $L$ can be extended to a linear function
from $M\otimes M$ to $R$, and this is an $R$-valued bicharacter of $M$.
This identifies the bicharacters of the group $L$ with the bicharacters of 
its group ring $M$. 

In order to define quantum vertex algebras we need a generalization of
commutative rings, called braided rings. The idea is that we should be
able to write $ab=\sum b_ia_i$ for suitable $a_i$ and $b_i$ related in
some way to $a$ and $b$. For example, for  a commutative ring we would
have $a_i=a$, $b_i=b$. The definition of the elements $a_i$ and $b_i$ is given 
in terms of an $R$-matrix with $R(a\otimes b)=\sum a_i\otimes b_i$, 
where an $R$-matrix is defined as follows. 
\proclaim Definition 2.3. An $R$-matrix for a ring $M$ 
with multiplication map $m:M\otimes M\mapsto M$
in an additive symmetric tensor category consists
of a map $R:M\otimes M\mapsto M\otimes M$ satisfying the following
conditions. 
\item{1.}
 ($R$ is compatible with 1.)
$R(1\otimes a)=1\otimes a$, $R(a\otimes 1)=a\otimes 1$. 
\item{2.}
  ($R$ is compatible with multiplication.) 
$m_{23} R_{12}R_{13} = R_{12}m_{23}:M\otimes M\otimes M\mapsto M\otimes M$ 
and
$m_{12} R_{23}R_{13} = R_{13}m_{12}:M\otimes M\otimes M\mapsto M\otimes M$.
\item{3.}
 (Yang-Baxter equation.) $R_{12}R_{13}R_{23}=R_{23}R_{13}R_{12}$.

Here $R_{13}$ is $R$ restricted to the first and third factors of
of $M\otimes M\otimes M$, and so on. 

\proclaim Definition 2.4. A braided ring $M$
in an additive symmetric tensor category  
is a ring $M$ with an $R$-matrix $R$ such that
$$  mR=m\tau:M\otimes M\mapsto M$$
(where $\tau:a\otimes b\mapsto b\otimes a$ is the twist map and
$m:a\otimes b\mapsto ab$ is the product). 

\noindent {\bf Example 2.5.}
Suppose that $M$ is $Z/2Z$ graded as
$M=M_0\oplus M_1$, and define $R$ by $R(a\otimes
b)=(-1)^{\deg(a)\deg(b)}a\otimes b$. Then $M$ is a braided ring with 
$R$-matrix $R$ if and only if $M$ is a super commutative ring. 

\proclaim Lemma/Definition 2.6. Suppose that $r$ is an $R$-valued
bicharacter
of a commutative cocommutative bialgebra $M$.
Define a new multiplication $\circ$ on $M$ by
$$a\circ b=\sum a'b'r(a''\otimes b'')$$
(where $\Delta(a)=\sum a'\otimes a''$, $\Delta(b)=\sum b'\otimes b''$). 
Then this makes $M$ into a  ring, called the twisting of
$M$ by $r$. If $r$ is symmetric 
then the twisting of $M$ by $r$ is
commutative. If $r$ is invertible (which is true whenever
$M$ is a Hopf algebra) then the twisting of $M$ by $r$ 
is a braided ring. 

Proof. The element 1 is an identity for twisting of $M$ by $r$
because $R$ is compatible with 1. The twisting is an associative
ring because $R$ satisfies the Yang-Baxter equation and is compatible
with multiplication.
It is easy to check that the twisting is commutative
if $r$ is symmetric and $M$ is commutative. 
Finally we have to check 
that $M$ has an $R$ matrix if $r$ is invertible. 
Define a bicharacter $r'$ by 
$$r'(a\otimes b)=\sum r(a'\otimes b')r^{-1}(b''\otimes a'').$$
We define the
$R$ matrix by
$$R(a\otimes b) = a'\otimes b'r'(b''\otimes a'')$$
where $\Delta(a) = \sum a'\otimes a''$,
and $\Delta(b) = \sum b'\otimes b''$.
It is easy to check that this satisfies the conditions
for an $R$ matrix for the twisting of $M$ by $r$. This proves lemma 2.6.

\noindent {\bf Example 2.7.}
Suppose that $L$ is a free abelian group or free abelian monoid with a
basis $\alpha_1,\ldots, \alpha_n$, and suppose that we are given
elements $r(\alpha_i,\alpha_j)\in S^*$ for some commutative
$R$-algebra $S$. We write $e^\alpha$ for the element
of the group ring of $L$ corresponding to the element
$\alpha\in L$.  Then we can extend $r$ to a unique $S$-valued
bicharacter of the  ring $M=R[L]$ by putting
$$r\big(\prod_i e^{m_i\alpha_i }\otimes  \prod_j e^{n_j\alpha_j}\big)
= \prod_{1\le i,j\le n} r(\alpha_i,\alpha_j)^{(m_i,n_j)}.
$$

\noindent {\bf Example 2.8.}
Suppose that $S$ is a commutative $R$ algebra and that $\Phi$ is a free
$R$-module, considered as an abelian Lie algebra.  We let $M$ be the
universal enveloping algebra of $\Phi$ (in other words the symmetric
algebra of $\Phi$), so $M$ is a commutative cocommutative Hopf algebra.  
Suppose that 
$r$ is any linear map from 
from $\Phi\otimes \Phi$ to $S$.
Then we can extend $r$ an $S$-valued bicharacter of  $M$ by putting
$$
r(\phi_1\cdots \phi_m\otimes \phi'_1\cdots \phi'_n)=
\cases
{
0& if $m\ne n$\cr
\sum_{\sigma\in S_m}\prod_{i=1}^m r(\phi_i\otimes 
\phi'_{\sigma(i)}) &if $m=n$\cr
}
$$
for $\phi_i,\phi'_i\in \Phi$ 
(where $S_m$ is the symmetric group of permutations of $1,2,\ldots,
m$).  

\noindent {\bf Example 2.9.}
Suppose that $r$ is any bicharacter of a cocommutative bialgebra
$M$. We can define an $R$-matrix for  $M$ by putting
$$R(a\otimes b) = \sum a'\otimes b' r(a''\otimes b'').$$
If $R$ is any $R$-matrix for a ring $M$ we can define a new associative
multiplication on $M$ by putting 
$$a\circ b=mR(a\otimes b):M\otimes M\mapsto M$$
where $m:M\otimes M\mapsto M$ is the old multiplication.
The composition of these two operations is just the twisting of $M$ by $r$. 

In the rest of this section we describe the construction of universal 
rings acted on by bialgebras, which we will need for the construction of
vertex algebras. These universal rings can be thought of as something like
the coordinate rings of function spaces or gauge groups. 
\proclaim Lemma/Definition 2.10. 
Suppose that $M$ is a commutative algebra over some ring and $H$ is a 
cocommutative coalgebra. 
Then there is a universal commutative algebra $H(M)$ 
such that there is a map $h\otimes m\mapsto h(m)$ from $H\otimes M$ to $H(M)$
with 
$$h(mn)=\sum h'(m)h''(n), \qquad h(1)=\eta(h).$$
If $H$ is a bialgebra then $H$ acts on the commutative ring
$H(M)$.  
If $M$ is a commutative and cocommutative bialgebra (or Hopf algebra)
then so is $H(M)$.

Proof. The existence of $H(M)$
is trivial; for example, we can construct it by writing down
generators and relations. 
Equivalently we can construct it as the quotient of the symmetric algebra
$S(H\otimes M)$ by the ideal generated by the images
of $H$ and $H\otimes M\otimes M$ under the maps describing the relations. 
If $H$ is a cocommutative bialgebra then it acts on $H(M)$ by 
$h_1(h_2(m))=(h_1h_2)m$. 
If $M$ has a coproduct $M\mapsto M\otimes M$
then this induces a map $M\mapsto H(M)\otimes H(M)$. As $H(M)\otimes H(M)$
is a commutative algebra acted on by $H$, this map extends to
a map from $H(M)$ to $H(M)\otimes H(M)$ by the universal property of $H(M)$.
It is easy to check that this coproduct makes $H(M)$ into a bialgebra.
This proves lemma 2.10. 

The ring $H(M)$ has the following geometric interpretation. 
Pretend that $H^*$ is the coordinate ring of a variety $G$. Then 
$Spec(H(M))$ can be thought  of as a sort of function space of all maps from 
$G$ to $Spec(M)$. If $H$ is a bialgebra then 
we can pretend that it is the group ring of a group $G$, and
the action of $H$ on $H(M)$ then corresponds to the
natural action of $G$ on this function space induced by 
the action of $G$ on itself by left multiplication. If in addition $M$
is a cocommutative Hopf algebra, then $Spec(M)$ is an affine algebraic group.
The space $Spec(H(M))$ is also an affine algebraic group, and can be thought
of as the gauge group of all maps from $G$ to $Spec(M)$. 

\noindent {\bf Example 2.11.}
Suppose that $H$ is the supercommutative bialgebra with a
basis $1$, $d$, with $d^2=0$, $\Delta(d)=d\otimes 1+1\otimes d$, such
that $d$ has odd degree. If $M$ is any supercommutative ring then
$H(M)$ is the ring of differential forms over $M$ (where of course
we replace ``commutative'' by ``supercommutative'' in lemma 2.10). 

\noindent {\bf Example 2.12.}
Suppose that $M$ is a polynomial algebra $R[\phi_1,\ldots,\phi_n]$.
Let $H$ be the commutative cocommutative Hopf algebra over $R$ with
basis $D^{(i)}$ for $i\ge 0$, where $D^{(i)}D^{(j)}={i+j\choose
i}D^{(i+j)}$ and $\Delta(D^{(i)})=\sum_jD^{(j)}\otimes D^{(i-j)}$.
(We can think of $H$ as the formal group ring of the one dimensional
additive formal group. If $R$ contains the rational numbers  
then $D^{(i)}=D^i/i!$ (where $D=D^{(1)}$) and $H$ is just the
universal enveloping algebra $R[D]$ of a one dimensional Lie algebra.)  Then
$H(M)$ is the ring of polynomials in the variables $D^{(i)}(\phi_j)$
for $i\ge 0$, $1\le j\le n$.
More generally, if we take $M$ to be a symmetric algebra $S^*(\Phi)$ for
an $R$-module $\Phi$, then $H(S^*(M))= S^*(H\otimes \Phi)$. 

\noindent {\bf Example 2.13.}
Suppose that $L$ is a lattice and $R[L]$ its group ring
and suppose that $H$ is the formal group ring of the one dimensional
additive group, as in example 2.12. Then $H(R[L])$ is the module
underlying the vertex algebra of the lattice $L$. 
If instead we take $H$ to be the polynomial ring $R[D]$
(with $\Delta(D)=D\otimes 1+1\otimes  D$) then  $H(R[L])$ is isomorphic to the 
tensor product $R(L)\otimes S^*(L(1)\oplus L(2)\oplus\cdots)$ of the
group ring $R[L]$ and the symmetric algebra of the sum
of an infinite number of copies $L(n)$ of $L\otimes
R$.  This tensor product is also commonly used to construct the
vertex algebra of a lattice. If $R$ contains the rational numbers then it
is equivalent to the first construction because $R[D]$ is then the
same as the $H$ defined in 2.12. However in non-zero characteristics 
it does not work quite so well; for example, we cannot define formal contour
integrals as in example 4.7, because this requires divided powers of $D$. 

We now show that bicharacters of $M$ are more or less the same
as $H\otimes H$-invariant bicharacters of $H(M)$. 

\proclaim Lemma 2.14. Suppose that $M$ and $N$ are bialgebras 
and $r$ is a bimultiplicative map from $M\otimes N$ to $S$, where $S$
is a commutative algebra acted on by the bialgebra $H$.  Then $r$
extends uniquely to a $H$ invariant bimultiplicative map from
$H(M)\otimes N$ to $S$. 

Proof. By adjointness we get an algebra homomorphism from $M$ to
the algebra $Hom(N, S)$ of linear maps from the coalgebra 
$N$ to the algebra $S$. By the universality property of $H(M)$ this extends
uniquely to a $H$ invariant homomorphism from $H(M)$ to $Hom(N,S)$,
which by adjointness gives a map from $H(M)\otimes N$ to $S$ such that
$r(m_1m_2\otimes n)=\sum r(m_1\otimes n')r(m_2\otimes n'')$ 
(where $\Delta(n)=\sum n'\otimes n''$).
To finish the proof we have to
check that $r(m\otimes n_1n_2)=\sum r(m'\otimes n_1)r(m''\otimes  n_2)$.  
The set of $m$
with this property contains $M$ because by assumption $r$ is
bimultiplicative on $M\otimes N$. It is also easy to check that
it is closed under
multiplication and under the action of $H$. Therefore it contains the
smallest $H$-invariant subalgebra of $H(M)$ containing $M$, which is the
whole of $H(M)$. This proves lemma 2.14.

\proclaim Lemma 2.15. 
Suppose that 
$H$ is a cocommutative bialgebra and 
$S$ is a commutative algebra acted on by $H\otimes H$.
Suppose that $M$ is a commutative and cocommutative 
bialgebra  with an $S$-valued bicharacter $r$. 
Then $r$ extends uniquely to a $H\otimes H$-invariant $S$-valued
bicharacter $r:H(M)\otimes H(M)\mapsto S$ of $H(M)$.

Proof. We apply lemma 2.14 to get a bimultiplicative $H$-invariant map
from $H(M)\otimes M$ to $S$. Then we apply lemma 2.14 again to get a
bimultiplicative $H\otimes H$-invariant map from $H(M)\otimes H(M)$ to
$S$. This proves lemma 2.15.

\proclaim 3.~Construction of some categories.

In this section we define a category $Fun(Fin^{\not\equiv},A, H, S)$ in
which we can carry out the ``twisted group ring'' 
construction in order to produce vertex algebras. 
The definition of this category is strongly motivated by
classical and quantum field theory, and commutative rings in this category
are formally quite similar to quantum field theories.

In the rest of this paper we fix
an additive tensor category $A$ that is cocomplete 
and such that colimits commute with tensor products. 
(In fact we do not need all colimits in $A$; it would be sufficient
for most applications to assume that $A$ has countable colimits.) For example, 
$A$ could be the category $R\hbox{-mod}$
of modules over a commutative ring. Note that most of the constructions
and definitions of section 2 work for any category $A$ with the properties 
above. 

\proclaim Definition 3.1. 
We define $Fin$ to be the category of all finite sets,
with morphisms given by functions. 
We define $Fin^{\not\equiv}$ to be the
category whose objects are finite sets with an equivalence relation $\equiv$, 
and whose morphisms are the functions $f$ preserving inequivalence;
in other words,  if $f(a)\equiv f(b)$
then $a\equiv b$. We define $\cup$ on $Fin$ and $Fin^{\not\equiv}$ 
to be the disjoint union (where in $Fin^{\not\equiv}$,  elements
of $I$ and $J$ in the disjoint union $I\cup J$ are inequivalent). 
This makes $Fin$ and $Fin^{\not\equiv}$ into (non-additive) symmetric 
tensor categories. 

We will write objects of $Fin^{\not\equiv}$ by using colons
to separate the equivalence classes. 

We could replace $Fin$ and $Fin^{\not\equiv}$ by smaller 
equivalent categories; for example we could restrict the objects
of $Fin$ to be the finite sets of the form $\{1,2,\ldots, n\}$. 

Note that $\cup $ is a coproduct in $Fin$ but is not a coproduct in
$Fin^{\not\equiv}$; in fact, $Fin^{\not\equiv}$ does not have
coproducts. For example the coproduct of a one point set and a two point set
with two equivalence classes does not exist. 

\proclaim Definition 3.2. If $C$ is a category we define the category
$Fun(C, A)$ to be the category of functors $V$ from $C$ to $A$.  
The category $Fun(C,A)$ is additive and has a 
symmetric tensor product given by the pointwise tensor product
$(U\otimes V)(I)=U(I)\otimes V(I)$. 

In applications the category $C$ will be one of $Fin$, $Fin^{\not\equiv}$,
or their opposite categories $Fin^{op}$, $Fin^{\not\equiv op}$ . 

\proclaim Definition 3.3. Suppose that $M$ is any commutative ring in $A$.  We
define $T_*(M)$ in $Fun(Fin, A)$ by $T_*(M)(I)=\otimes _{i\in I}M$, where
the action of $T_*(M)$ on morphisms of $Fin$ is induced in the obvious
way by the product and unit of $M$.  (For example, if $f$ is the
morphism from $\{1,2\}$ to itself with $f(1)=f(2)=2$, then $T_*(M)(f)$
takes $x_1\otimes x_2$ to $1\otimes x_1x_2$.)  If $H$ is a
cocommutative coalgebra in $A$ 
then we define $T^*(H)$ in $Fun(Fin^{op}, A)$ in a
similar way, using the coproduct and counit of $H$ to define the action of
$T^*(H)$ on morphisms.

\noindent {\bf Example 3.4.}
If $M$ is a commutative ring in $A$ then $T_*(M)$ is a commutative ring
in $Fun(Fin, A)$. If in addition $M$ is a commutative cocommutative
bialgebra then so is $T_*(M)$.

\noindent {\bf Example 3.5.}
If $V$ is a commutative ring in $Fun(Fin^{\not\equiv}, A)$ then 
we can  think of $V(I)$ as the space of (nonsingular)
quantum fields $\phi(x_1,x_2,\ldots)$
depending on $|I|$ spacetime variables.

The space of fields in one spacetime variable is acted on by the group
of automorphisms $G$ of spacetime, and similarly the space of fields
of $|I|$ spacetime variables is acted on by $|I|$ commuting copies of
$G$. We now add a similar structure to the objects of $Fun(Fin,
A)$. It is convenient to use a cocommutative bialgebra $H$ instead
of a group $G$; we can think of this bialgebra $H$ as analogous to
the group ring of the automorphisms of spacetime (or maybe to the
universal enveloping algebra of the Lie algebra of infinitesimal
automorphisms of spacetime).

\proclaim Definition 3.6. Suppose that $T$ is a cocommutative bialgebra
in $Fun(C^{op}, A)$. (In applications, $T$ will be of the form $T^*(H)$
for a cocommutative bialgebra of $A$.)
We define a
$T$ module in $Fun(C,A)$ to be an object $V$ of $Fun(C,A)$ such
that $V(I)$ is a module over $T(I)$ for all $I$ and such that
$f_*(f^*(g)(v))= g(f_*(v))$ for $v\in V(I)$, $g\in T(J)$,
$f:I\mapsto J$. The action of $T$ on the tensor product of two
$T$ modules is defined in the usual way
using the coalgebra structure of the $T(I)$'s.
We define $Fun(C, A, T)$ to be the additive symmetric tensor category
of $T$ modules in $Fun(C, A)$.

\noindent {\bf Example 3.7.}
Suppose that $V$ is any commutative ring in $A$ acted on by the
cocommutative bialgebra $H$. Then  $T_*(V)$  is a commutative
ring in $Fun(Fin,A,T^*(H))$. 

Recall that we can define the category of modules over any commutative
ring in any additive symmetric tensor  category, and it is again 
an additive  symmetric tensor category. 

\proclaim Definition 3.8. Suppose  that $T$ is a cocommutative bialgebra
in $Fun(C^{op},A)$
and suppose that $S$ is a commutative ring in $Fun(C, A,T)$. We define
$Fun(C,A,T,S)$ to be the additive symmetric tensor category of modules
over $S$.

\noindent {\bf Example 3.9.}
Suppose that we define $S$ by letting $S(I)$ be the smooth
functions depending on $|I|$ variables in spacetime. Then we would
expect a field theory to be a module over $S$ because we should be
able to multiply a field by a smooth function to get a new field.

Commutative rings in $Fun(Fin,A, T^*(H), S)$ as defined above 
behave rather like classical field theories, or at least they have
most of their formal properties. However quantum field theories
do not fit into this framework. The problem is that
in quantum field theory it is no longer true that the product of
two nonsingular fields is a nonsingular field. For example, a typical formula
in free quantum field theory is
$$\phi(x_1)\phi(x_2)=:\phi(x_1)\phi(x_2):+\Delta(x_1-x_2)$$ 
where the propagator $\Delta(x)$
usually has a singularity at $x=0$. In particular if we take $x_1=x_2$ we
find that the product of two fields depending on $x_1$ is not
defined. Instead, we can take the product of two fields
depending on different variables $x_1$ and $x_2$, and it lies
in the space $V(1:2)$ of fields that are defined whenever
$x_1$ and $x_2$ are ``apart'' in some sense. 

The category $Fun(C,A,T^*(H), S)$ has a natural tensor product
$\otimes$ which can be used to define multilinear maps. We will now
define a new tensor product in $Fun(C,A,T^*(H),S)$ by
defining a new concept of multilinear maps, called singular multilinear maps.
We assume that $C$ is a symmetric tensor category (not necessarily additive)
with the tensor product denoted by $\cup$. (As the notation suggests, 
this will often be some sort of disjoint union.) 

\proclaim Definition 3.10. 
We let $T$ be a cocommutative bialgebra in $Fun(C^{op},A)$, and we let
$S$ be a commutative ring in $Fun(C,A,T)$. 
Suppose that $U_1,U_2,\ldots$ and $V$
are objects of $Fun(C,A,T,S)$. We define a singular multilinear map from
$U_1,U_2,\ldots $ to $V$ to be a set of maps from $U_1(I_1)\otimes_A
U_2(I_2)\cdots$ to $V(I_1\cup I_2\cdots)$ for all
$I_1,I_2,\ldots \in C$, satisfying the following
conditions.
\item{1} The maps commute with the action of $T$.
\item{2} The maps commute with the actions of $S(I_1)$, $S(I_2),\ldots$. 
\item{3} If we are given any morphisms from $I_1$ to $I_1'$, 
$I_2$ to $I_2', \ldots$, then the following diagram commutes:
$$\matrix{
U_1(I_1)&\otimes& U_2(I_2)  \cdots 
&\longrightarrow 
& V(I_1\cup I_2\cdots)\cr
&\downarrow& &&\downarrow \cr
U_1(I_1')&\otimes& U_2(I_2')  \cdots 
&\longrightarrow 
& V(I_1'\cup I_2'\cdots)\cr
}$$

As $A$ is co-complete and co-limits commute with taking tensor products
the singular multilinear maps are representable, so we  define the
``singular tensor products'' $U_1\odot U_2\cdots$ to be the objects
representing the singular multilinear maps. It is possible 
to write down an explicit formula for
these singular tensor products as follows.
$$(U_1\odot U_2\cdots)(I)
= \lim_{\longrightarrow\atop I_1\cup I_2\cdots\mapsto I} 
\left(U_1(I_1)\otimes_AU_2(I_2)\cdots\right)
\otimes_{S(I_1)\otimes S(I_2)\cdots}S(I)$$
where the limit is a direct limit taken over the following category.
The objects $I_1\cup I_2\cdots\mapsto I$ of the category consist of objects $I_1,I_2,\ldots$ of
$C$ together with a morphism from $I_1\cup I_2\cdots$ to $I$.
A morphism from $I_1\cup I_2\cdots \mapsto I$ to $I_1'\cup I_2'\cdots\mapsto I$
consists of morphisms from $I_1$ to $I_1'$, $I_2$ to $I_2'$, $\ldots$, 
making the following diagram commute:
$$\matrix{
I_1\cup I_2&\cdots &\longrightarrow& I\cr
\downarrow&&&\Vert\cr
I_1'\cup I_2'&\cdots &\longrightarrow& I.\cr
}$$

J. M. E. Hyland told me that the product $\odot$ is similar to the
``Day product'' in category theory.  The construction of $\odot$ can
be extended to the case when $C$ is a ``symmetric multi-category''
rather than a symmetric tensor category.  Soibelman remarked that the
conditions for $V$ to be an algebra for $\odot$ are similar to the
conditions for the functor $V$ from $C$ to $A$ to be a functor of
tensor categories.

\noindent {\bf Example 3.11.}
Suppose that $\cup $ is a coproduct in $C$; for example, 
we could take $C$ to be $Fin$ and $\cup$ to be disjoint union. 
Then singular tensor products are the same as pointwise
tensor products.  In later examples we will take $C$ to be
$Fin^{\not\equiv}$ and $\cup$ to be disjoint union, which is not a
coproduct in $Fin^{\not\equiv}$.

The two tensor products $\odot$ and $\otimes$ are related in several ways,
as follows.  
There is a canonical morphism from $U\odot V$ to $U\otimes V$, 
so that any ring is automatically a singular ring. 
Also there is a canonical ``interchange'' morphism
$$(U\otimes V)\odot (W\otimes X)\rightarrow (U\odot W)\otimes (V\odot
X).$$ (Unlike the case of the interchange map for natural
transformations, this interchange map is not usually an
isomorphism.)  The interchange map can be used to show
that if $U$ and $V$ are singular rings then so is $U\otimes V$.

We define singular rings, singular Lie algebras, and so on, in 
$Fun(C, A, T, S)$
to be rings, Lie algebras, and so on using the singular tensor product. 
We define singular bialgebras a little bit differently: the product
uses the singular tensor product, but the coproduct uses the pointwise
tensor product $\otimes$. Note that for this to make sense we need to know that
the pointwise tensor product of two singular algebras is a singular algebra; 
see the paragraph above. In general, we should use the pointwise
tensor product $\otimes$ for ``coalgebra'' structures, and
the singular tensor product $\odot$ for ``algebra'' structures.

If $S$ is a commutative ring in $Fun(Fin^{\not\equiv}, A, T^*(H))$ then by
restriction it is also a commutative ring in $Fun(Fin, A, T^*(H))$ (using
the functor which gives any finite set the equivalence relation where
all elements are equivalent.) We can embed the category $Fun(Fin, A,
T^*(H), S)$ into $Fun(Fin^{\not\equiv}, A, T^*(H), S)$ by defining
$V(I_1:I_2\cdots)=V(I_1\cup I_2\cdots)\otimes_{S(I_1)\otimes
S(I_2)\cdots}S(I_1:I_2\cdots)$ for $I_1,I_2,\ldots\in Fin$. In particular singular multilinear
maps are defined in $Fun(Fin, A, T^*(H), S)$.
(Note that singular tensor products representing singular
multilinear maps do not usually exist in  $Fun(Fin, A, T^*(H), S)$,
though they do exist in the larger category $Fun(Fin^{\not\equiv}, A, T^*(H),
S)$.)

The main point of all this category theory is the following definition:
\proclaim Definition 3.12. 
Suppose that $A$ is an additive symmetric tensor category, $H$ is a
cocommutative bialgebra in $A$, and $S$ is a commutative ring in
$Fun(Fin^{\not\equiv}, A, T^*(H))$.  We define an $(A, H, S)$ vertex
algebra to be a singular commutative ring in $Fun(Fin, A, T^*(H), S)$.
We define a quantum $(A, H,S)$ vertex algebra to be a singular braided
ring in $Fun(Fin,A, T^*(H), S)$.

Soibelman remarked that all the examples of quantum $(A, H,S)$ vertex
algebras in this paper have the extra property that the
$R$ matrix satisfies $R_{12}R_{21}=1$, so perhaps this condition
should be added to the definition of a quantum $(A, H,S)$ vertex algebra. 

Note that the vertex algebra is in $Fun(Fin, A, T^*(H), S)$
rather than  $Fun(Fin^{\not\equiv}, A, T^*(H), S)$, although we can of
course embed the former category in the latter if we wish. 
The reason for using 
$Fun(Fin, A, T^*(H), S)$ rather than $Fun(Fin^{\not\equiv}, A, T^*(H), S)$
is that we wish to have control over the connection between (say)
$V(1,2)$ and $V(1:2)$. 

\proclaim  4.~Examples of  vertex algebras. 

In this section we construct some examples of $(A,H,S)$ vertex algebras
by applying the twisted group ring construction of section 2 to 
the categories constructed in section 3. We also show how these
are related to classical vertex algebras. 

\proclaim Lemma 4.1. Suppose that $r$ is an $H\otimes H$-invariant 
$S(1:2)$-valued bicharacter 
of a commutative cocommutative bialgebra $H(M)$ in $A$. Then $H$ 
can be extended to a singular bicharacter of $T_*(H(M))$, which we also denote
by $r$. 

Proof. We define $r$ by
$$r(\bigotimes_{i\in I}a_i\otimes\bigotimes_{j\in J}b_j)
= \sum \prod_{i\in I}\prod_{j\in J} r(a_i^{(j)}\otimes b_j^{(i)})$$
where 
$\Delta^{|J|-1}(a_i)=\sum \bigotimes_{j\in J}a_i^{(j)}$,
$\Delta^{|I|-1}(b_j)=\sum \bigotimes_{i\in I}b_j^{(i)}$,
and $r(a_i^{(j)}\otimes b_j^{(i)})$ is considered as
an element of $S(I\cup J)$ using the obvious map from
$S(i:j)$ to $S(I\cup J)$. Some routine checking then proves lemma 4.1.

The following theorem is the main theorem of this paper. It shows how
to construct many examples of $(A, H, S)$ vertex algebras, by
giving a sort of generalization of the construction of the vertex algebra of
a lattice. 

\proclaim Theorem 4.2. Suppose that $H$ is a cocommutative bialgebra in $A$
and $S$ is a commutative ring in $Fun(Fin^{\not\equiv}, A, T^*(H))$.
Assume that we are given an $S(1:2)$-valued bicharacter $r$ of a
commutative and cocommutative bialgebra $M$ in $A$.  The bicharacter
$r$ of $M$ extends to a bicharacter of $T_*(H(M))$ as in lemmas 2.15
and 4.1, which we also denote by $r$. Then the twisting of $T_*(H(M))$
by $r$ is a quantum $(A,H,S)$ vertex algebra if $r$ is invertible, and
is an $(A,H,S)$ vertex algebra if $r$ is symmetric.

Proof. By lemma 2.10  and example 3.4, $T_*(H(M))$ is a commutative cocommutative
bialgebra in $Fun(Fin, A, T^*(H), S)$. By lemmas 2.15 and 4.1
 the bicharacter $r$ extends to a singular bicharacter of $T_*(H(M))$ with
values in $S$. By lemma 2.6 (extended to additive tensor categories)
the twisting of $T_*(H(M))$ by $r$ is a braided ring if $r$ is invertible, 
and  is
a commutative ring if $r$ is symmetric.  Theorem 4.2 now follows from the
definition 3.12 of (quantum) $(A,H,S)$ vertex algebras.

The following theorem describes the relation between the
$(A, H, S)$ vertex algebras of this paper, and ordinary vertex algebras. 

\proclaim Theorem 4.3. 
Suppose we take $H$ to be the formal group ring of the one dimensional
additive formal group, as in example 2.12. Define $S$ by $S(I)=$ the
$R$-algebra generated by $(x_i-x_j)^{\pm 1}$ for $i$ and $j$ not equivalent
(so $S=R$ if all elements of $I$ are equivalent).  
If $V$ is a $(R\hbox{-mod}, H, S)$ vertex algebra, then $V(1)$ is an ordinary
vertex algebra over the ring $R$. 

Proof. For every element $u_1$ of $V(1)$ we have to
construct a vertex operator $u_1(x_1)$ taking $V(1)$ to
$V(1)[[x_1]][x_1^{-1}]$. We do this as follows.  If $u_2\in V(2)$ then
$u_1u_2\in V(1:2)=V(1,2)\otimes S(1:2)=V(1,2)[(x_1-x_2)^{\pm 1}]$.
There is a map from $V(1,2)$ to $V(1)[[x_1, x_2]]$ taking $w$ to the
``Taylor series expansion'' $\sum_{i,j}f_{12\rightarrow
1}(D_1^{(i)}D_2^{(j)}w)x_1^ix_2^j$. (Here $f_{12\rightarrow 1}$ is the
map from $V(1,2)$ to $V(1)$ induced by the morphism of finite sets
taking both 1 and 2 to 1, and $D_1$ and $D_2$ indicate the two 
different actions of $H$ on $V(1,2)$.) This induces a map from 
$V(1,2)[(x_1-x_2)^{\pm 1}]$ to $V(1)[[x_1,x_2]][(x_1-x_2)^{-1}]$, 
and we denote the image of $u_1u_2 $ under this map by $u_1(x_1)u_2(x_2)$. 
Then we define the vertex operator $u_1(x_1)$ by
$u_1(x_1)u_2= u_1(x_1)u_2(0)\in V(1)[[x_1]][x_1^{-1}]$. 

This defines the vertex operators of elements of $V(1)$; now we have
to check that they formally commute. 
We can define expressions like $u_1(x_1)u_2(x_2)u_3(x_3)\cdots
\in V(1)[[x_1,\ldots]][\prod(x_i-x_j)^{-1}]$ in the same way as above.
The fact that $V$ is commutative implies that 
$u_1(x_1)u_2(x_2)u_3(0)=u_2(x_2)u_1(x_1)u_3(0)$. This in turn implies
that the vertex operators $u_1(x_1)$ and $u_2(x_2)$ commute in the sense
that 
$$(x_1-x_2)^N (u_1(x_1)u_2(x_2)-u_2(x_2)u_1(x_1))u_3=0$$
for $N$ a sufficiently large integer, depending on $u_1$ and $u_2$. 
So we have constructed commuting vertex operators for all elements
of $V(1)$, and this can easily be used to show that
$V(1)$ is a vertex algebra. This proves theorem 4.3.

\noindent {\bf Example 4.4.}
Take $L$ to be an even integral lattice. Choose a bicharacter
$c$ such that $c(\alpha, \beta)=(-1)^{(\alpha,\beta)}c(\beta,\alpha)$.
(There are many ways to do this. For example we 
can choose a basis $\alpha_1,\alpha_2,\ldots$ and define
$c(\alpha_i,\alpha_j)$ to be 1 if $i\ge j$ and $(-1)^{(\alpha_i,\alpha_j)}$
if $i<j$.) 
Define a symmetric $R[(x_1-x_2)^{\pm1}]$-valued bicharacter $r$ of $L$
by
$$r(\alpha,\beta)(x_1,x_2)=c(\alpha,\beta)(x_1-x_2)^{(\alpha,\beta)}.$$
If $V$ is the $(R\hbox{-mod}, H, S)$ vertex algebra constructed in theorem 4.2
with underlying object $T_*(H(R[L]))$ then $V(1)$ is just the usual
vertex algebra of the even integral lattice $L$.  If $L$ is any
integral lattice (not necessarily even) then we can do a similar
construction with the following changes.  We choose $c$ so that
$c(\alpha,
\beta)=(-1)^{(\alpha,\beta)+(\alpha,\alpha)(\beta,\beta)}c(\beta,\alpha)$.
The bicharacter $r$ is no longer symmetric but is supersymmetric,
so we end up with a vertex superalgebra rather than a vertex algebra.

\noindent {\bf Example 4.5.}
Now we write down some quantum deformations of example 4.4.
Let $L$ be an even lattice as in example  4.4, 
let $q$ be an invertible element of the commutative  ring $R$, 
and let $A$ be the category of
$R$ modules. 
We define $S$ by $S(I)=R$ if $I$ has only one equivalence class, 
and $S(I)=$ the $R$-algebra generated by
$(x_i-q^nx_j)$ for $i$ and $j$ not equivalent, $n$ an integer,
if $I$ has more than 1 equivalence class. 
Choose a basis
$\alpha_1\ldots,\alpha_n$ for $L$ and define $r$ using lemma 2.7
 by putting
$$
r(\alpha_i,\alpha_j)=
c(\alpha,\beta)
\prod_{k=1}^{(\alpha_i,\alpha_j)} (x_1-q^{(\alpha_i,\alpha_j)-2k}x_2)
$$
where $c$ is the bicharacter of example 4.4. 
By applying theorem 4.2 we get a $(R\hbox{-mod}, H, S)$ quantum vertex algebra. 
We see that
$$
(x_1-q^{(\alpha_i,\alpha_j)}x_2)
e^{\alpha_1 }(x_1)e^{\alpha_2 }(x_2) 
= (q^{(\alpha_i,\alpha_j)}x_1-x_2)
e^{\alpha_2 }(x_2)e^{\alpha_1 }(x_1).
$$ 
This is similar to many of the formulas of statistical mechanics
in the book [J-M].

\noindent {\bf Example 4.6.}
We show how to construct $(A,H,S)$ vertex algebras corresponding
to generalized free quantum field theories. 
Suppose that $\Phi$ is a module over a commutative ring
in $A$ and $H$ is a commutative cocommutative bialgebra in 
$A$. Then any linear map $\Delta$ from $\Phi\otimes\Phi$ to $S(1:2)$ gives a
quantum $(A,H,S)$ vertex algebra as follows.  Use example 2.8  to extend $r$ to
a $S(1:2)$-valued bicharacter of the symmetric algebra $M$ of $\Phi$.
Then use theorem 4.2 to make
$T_*(H(M))$ into a quantum $(A,H,S)$ vertex algebra. 
If $r$ is symmetric then this
is a $(A,H,S)$ vertex algebra, and is  closely related to generalized
free quantum field theories, at least when $H$ is finite dimensional
abelian. (To obtain analogues of free quantum field theories in odd
dimensions or dimension 2 we should allow slightly more general sorts
of singularities, such as half integral powers or logarithms of
$(x_1-x_2)^2$ rather than just poles.)  The function $r$ gives the
propagator of free fields, and the Greens functions 
$\langle|\phi_1(x_1)\cdots\phi_n(x_n)|\rangle$ can be recovered
as $\eta(\phi_1(x_1)\cdots\phi_n(x_n))$ where $\eta$ is the counit of
$H(M)$ and $\phi_1,\ldots,\phi_n $ are elements of $\Phi$. 

Take $H$ to be the additive formal group of dimension $d$ for some
positive even integer $d$.  If we take $\Phi$ to be a one dimensional
free module over $R$ spanned by an element $\phi$ and put
$r(\phi\otimes \phi) = (\sum(x_{1,i}-x_{2,i})^2)^{1-d/2}$ then $V$ is the ``$H$
vertex algebra of a free scalar field'' constructed in [B98].  It is obvious
that we can just write down many quantum deformation of this $H$
vertex algebra just by varying $r$; for example, we could take
$r(\phi\otimes \phi) = (\sum(x_{1,i}-qx_{2,i})^2)^{1-d/2}$.

\noindent {\bf Example 4.7.}
In the theory of vertex algebras we often get contour integrals
such as 
$$\int_{x_1}a_1(x_1)a_2(x_2)a_3(x_3)dx_1.$$
We will show how to define such contour integrals for
$(A, H, S)$ vertex algebras, where $H$ and $S$ are as in theorem 4.3. 
Take $a_i\in V(i)$, where $V(i)$ can be identified with $V(1)$. 
We know that $a_1a_2a_3\in V(1:2:3)$ using the multiplication of
$V$. We also know that $V(1:2:3)=V(1,2,3)[(x_1-x_2)^{\pm 1},
(x_2-x_3)^{\pm 1},(x_1-x_3)^{\pm 1}]$, so we can write
$a_1a_2a_3$ as a finite sum of terms of the form
$$a_{123}(x_1-x_2)^i(x_1-x_3)^j(x_2-x_3)^k.$$
Next we expand $(x_1-x_3)^j$ as a possibly infinite
series
$$(x_1-x_3)^j=\sum _{n\ge 0} {j\choose n} (x_1-x_2)^n(x_2-x_3)^{j-n}.$$
Finally we replace each term $a_{123}(x_1-x_2)^i(x_2-x_3)^k$
by
$$f_*(D_1^{(i)}(a_{123}))(x_2-x_3)^k\in V(2:3)$$ where $f$ is the
function from $\{1,2,3\}$ to $\{2,3\}$ with $f(1)=f(2)=2$, $f(3)=3$.
This algebraically defined contour integral has most of the properties
one would expect. For example we have the identity
$$
\eqalign{
&\int a_1(x_1)dx_1\int a_2(x_2)dx_2a_3-\int a_2(x_2)dx_2\int a_1(x_1)dx_1a_3\cr
=& \int \left(\int a_1(x_1)dx_1a_2\right)(x_2)dx_2 a_3\cr
}
$$
which can be used to prove  the usual vertex algebra identities. 
Of course this identity depends on  the simple choice of $H$ and $S$ we made;
for more complicated choices of $H$ and $S$ we will usually 
get more complicated identities.
In particular contour integrals can be defined in terms of the more
elementary operations of a $(A,H,S)$ vertex algebra.  
One reason for using the bialgebra $H$ with divided powers (see example
2.12) rather then the universal enveloping algebra $R[D]$ is that the
divided powers are needed to define the contour integrals. 

\noindent {\bf Example 4.8.}
Take $H$ as in example 2.12, and let $S(I)$ be $R$ if $I$ has at most
one equivalence class, and the ring generated by the
elements $(x_i-q^nx_j)^{\pm 1}$ for $i\not\equiv j$ and $I$ having more than 
1 equivalence class.  Then if $V$ is a
quantum $(A, H, S) $ vertex algebra, we can think of $V(1)$ as being
some sort of ``quantum vertex algebra''.  We will not give a
definition of quantum vertex algebras here, because the philosophy of
this paper is that (quantum) vertex algebras should be replaced by
(quantum) $(A,H,S)$ vertex algebras.  Several sets of axioms for
quantum vertex algebras have been proposed by various authors in
[E-K], [F-R].

\noindent {\bf Example 4.9.}
The (ordinary) tensor product of two (ordinary) vertex algebras is a
vertex algebra.  The analogue of this for $(A,H,S)$ vertex algebras is
trivial to prove: the pointwise tensor product of any two singular commutative
rings in $Fun(Fin^{\not\equiv}, A, H, S)$ is a singular commutative
ring, and the pointwise tensor product of two objects of $Fun(Fin, A,
T^*(H), S)$ is still in $Fun(Fin, A, T^*(H),S)$, so the pointwise
tensor product
of two $(A, H,S)$ vertex algebras is an $(A, H,S)$ vertex algebra.
Note that the singular tensor product of two $(A,H,S)$ vertex algebras
is  a singular commutative ring in $Fun(Fin^{\not\equiv}, A,
T^*(H), S)$, but need not be in $Fun(Fin, A, T^*(H),S)$, so the
singular tensor product of two $(A,H,S)$ vertex algebras need not be
an $(A,H,S)$ vertex algebra.

\noindent {\bf Example 4.10.}
We can obtain many variations of vertex algebras by changing
$H$ and $S$. For example we could take $H$ to be the 
universal enveloping algebra of the Virasoro algebra to 
get things similar to ``vertex operator algebras''. 
If we take $H$ to be the tensor product of two copies of the Virasoro 
algebra then $(A,H,S)$ vertex algebras are closely related to
conformal field theory and string theory. If we let $H$ 
be the universal enveloping algebra of various  superalgebras
then we get $(A,H,S)$ vertex algebras related to supersymmetry.

\proclaim 5.~Open problems.

In this section we  list some suggestions for further research.

\noindent {\bf Problem 5.1.} 
Are there natural quantum deformations of other well known vertex
algebras, such as the monster vertex algebra [B86], [F-L-M], the
vertex algebra of the lattice $II_{25,1}$ [B86], [K97], and the vertex
algebras of highest weight representations of affine Lie algebras and
the Virasoro algebra [F-Z], [K97]?  Etingof and Kazhdan [E-K]
construct ``quantum vertex operator algebras'' corresponding to the
vertex algebras of affine Lie algebras, and it seems likely that their
construction could be extended to give examples satisfying the
definitions in this paper.  Frenkel and Jing [F-J] previously
constructed vertex operators related to of quantum affine Lie
algebras.

\noindent {\bf Problem 5.2.} 
Ordinary vertex algebras can be used to construct many examples of
generalized Kac-Moody algebras. Is there a relation between quantum
vertex algebras and some sort of quantized generalized Kac-Moody
algebras, possibly those defined in [K95]?

\noindent {\bf Problem 5.3.} 
The similarity of the formulas in solvable lattice models in [J-M] and
quantum vertex algebras suggests that there may be some relation
between these subjects.

\noindent {\bf Problem 5.4.} 
We have constructed vertex algebras from bicharacters of bialgebras
that are both commutative and cocommutative. If a bialgebra is
cocommutative but not commutative then the bicharacters are usually
not all that interesting (for the much same reason that one
dimensional characters of a non-abelian group are not usually
interesting). However there are nontrivial examples of bicharacters of
bialgebras that are neither commutative or cocommutative. Can these
be used to construct some sort of vertex algebras?

\noindent {\bf Problem 5.5.} 
Construct $(R\hbox{-mod}, H, S)$ vertex algebras corresponding
to the other standard examples of vertex algebras, such as
the vertex algebras of affine and Virasoro algebras ([F-Z]), or
the monster vertex algebra ([F-L-M])  or the vertex algebra of
differential operators on a circle ([K97]). 

\noindent {\bf Problem 5.6.} 
Many of the constructions and definitions in section 3 do
not use the fact that the category $A$ is additive. Is there
any use for these constructions in the non-additive case? 

\noindent {\bf Problem 5.7.} 
Do these constructions for braided rather than symmetric
tensor categories. In particular it should be possible
to allow nonintegral powers of $x_i-x_j$, which often arise
from non-integral lattices or from conformal field theory. 

\noindent {\bf Problem 5.8}
A cobraided Hopf algebra (as defined in in [K, definition
VIII.5.1]) is a
Hopf algebra with a bicharacter $r$ with the extra property that
$\mu^{op}=r*\mu*\bar r$. This suggests that it might be possible
to replace commutative, cocommutative bialgebras by
something more general, maybe cobraided bialgebras. 
In particular theorem 4.2 should be extended to the case when 
$M$ is cobraided rather than cocommutative. 

\noindent {\bf Problem 5.9} 
Instead of twisting a group ring by a bicharacter,
we can also twist it by a 2-cocycle (preferably normalized). 
We can define ``multiplicative 2-cocycles'' of arbitrary
cocommutative bialgebras with values in any algebra $S$ acted on by
the bialgebra, and use these to construct more general twistings.
We can also define multiplicative $n$-cochains, cocycles, and
coboundaries, and use these to define multiplicative analogues
$H^n(M,S^*)$ of cohomology groups.  Note that the usual (additive)
cohomology $H^n(M,S)$ of bialgebras depends only on the underlying
associative algebra and the counit of $M$ and on the module structure
of $S$, and should not be confused with these multiplicative
cohomology groups $H^n(M, S^*)$ that also depend on the coproduct of
$M$ and the algebra structure of $S$.
Find some examples of vertex algebras constructed using singular 2-cocycles
rather than singular bicharacters. There are many examples 
that can be constructed like this in a formal (and not very interesting) way
from a perturbative quantum field theory.

\noindent {\bf Problem 5.10} 
It is possible to construct singular 2-cocycles
which look formally similar to the Greens functions of perturbative
quantum field theories. At the moment this just seems to be little
more than a formal triviality, but may be worth investigating further.

\noindent {\bf Problem 5.11} 
I. Grojnowski and S. Bloch independently suggested replacing the Hopf
algebra $H$ of example 4.4 by the formal group ring of the formal
group of an elliptic curve. Over the rationals this makes no
difference, but over finite fields or the integers
we seem to get something different.
The underlying space of the vertex algebra we get can be thought of as
the coordinate ring of the gauge group of 
maps from a (formal) elliptic curve to an
algebraic torus. The problem is to find a use for this construction!

\noindent {\bf Problem 5.12} 
Develop the theory of categories with two symmetric tensor products
satisfying the conditions suggested in section 3
(and maybe some others), and find more examples of
them. Soibelman pointed out that Beilinson and Drinfeld [B-D] have 
some categories which have both a tensor product and a separate multilinear
structure.

\noindent {\bf Problem 5.13} 
The study of orbifolds of vertex algebras (in other words, fixed
subalgebras under finite automorphism groups) is notoriously hard (see
[D-M] for example), though this ought to be an easy and natural
operation.  The difficulties appear to be caused partly by the fact
that vertex algebras seem to have something missing from their
structure.  Does the theory of orbifolds for $(A,H,S)$ vertex algebras
(with their extra structure of fields of several spacetime variables)
become any easier?

\noindent {\bf Problem 5.14}
Soibelman suggested that the examples of associative algebras of
automorphic forms in the meromorphic tensor category of [So, Theorem
8] might be some sort of $(A,H,S)$ vertex algebras. These may be
related to the algebras in [K96].

\proclaim References.

\item{[B86]} R. E. Borcherds,
Vertex algebras, Kac-Moody algebras and the monster, 
Proc.Nat. Acad. Sci. U.S.A. 83 (1986), 3068-3071.
\item{[B98]} R. E. Borcherds,
Vertex algebras.  in ``Topological field theory, primitive forms and
related topics'' (Kyoto, 1996), edited by M. Kashiwara, A. Matsuo,
K. Saito and I. Satake.  Progress Math., 160, Birkh\"auser Boston,
Inc., Boston, MA. ISBN: 0-8176-3975-6 pp. 35--77.  Birkh\"auser Boston,
Boston, MA, 1998.  Preprint q-alg/9706008.
\item{[B-D]} A. Beilinson, V. Drinfeld, Chiral algebras I. 1995(?) preprint.
\item{[D-M]} C. Dong, G.  Mason, 
On quantum Galois theory. 
Duke Math. J. 86 (1997), no. 2, 305--321. 
\item{[E-K]} P. Etingof, D. Kazhdan, Quantization of Lie bialgebras, V.
Preprint math.QA/9808121
\item{[F-R]} E. Frenkel, N. Yu. Reshetikhin, Towards Deformed Chiral Algebras,
preprint \hbox{q-alg/9706023}, submitted to the Proceedings of the Quantum
Group Symposium at the XXIth International Colloquium on Group
Theoretical Methods in Physics, Goslar 1996.
\item{[F-J]} I. B. Frenkel, N. H.  Jing,
Vertex representations of quantum affine algebras. Proc. Nat.
Acad. Sci. U.S.A. 85 (1988), no. 24, 9373--9377.
\item{[F-L-M]}
I. B. Frenkel, J.  Lepowsky, A.  Meurman, Vertex operator algebras and
the Monster. Pure and Applied Mathematics, 134. Academic Press, Inc.,
Boston, MA, 1988. ISBN: 0-12-267065-5
\item{[F-Z]} I. B. Frenkel, Y. Zhu, 
Vertex operator algebras associated to representations of
affine and Virasoro algebras. Duke Math. J. 66 (1992), no. 1, 123--168. 
\item{[J-M]} M. Jimbo, T.  Miwa, 
Algebraic analysis of solvable lattice models. CBMS Regional
Conference Series in Mathematics, 85. 
Published for the Conference Board of the Mathematical Sciences,
Washington, DC; by the American Mathematical Society, Providence, RI, 1995. 
ISBN: 0-8218-0320-4 
\item{[J]} N. Jing, 
Quantum Z-algebras and representations of quantum affine algebras, 
preprint math.QA/9806043.
\item{[K97]} V. Kac, 
Vertex algebras for beginners. 
University Lecture Series, 10. American Mathematical
Society, Providence, RI, 1997.  ISBN: 0-8218-0643-2
\item{[K95]} S.-J. Kang,
Quantum deformations of generalized Kac-Moody algebras and their modules.
J. Algebra 175 (1995), no. 3, 1041--1066. 
\item{[K96]} M. M. Kapranov, 
Eisenstein series and quantum affine algebras, preprint alg-geom/9604018.
\item{[K]} C. Kassel, 
Quantum groups. Graduate Texts in Mathematics, 155. Springer-Verlag, New
York, 1995.  ISBN: 0-387-94370-6 
\item{[L]} J. Lambek, 
Deductive systems and categories. II. Standard constructions and
closed categories. 1969 Category Theory, Homology Theory and their
Applications, I (Battelle Institute Conference, Seattle, Wash., 1968,
Vol. One) pp. 76--122 Springer, Berlin
\item{[So]} Y. Soibelman, 
Meromorphic tensor categories, preprint q-alg/9709030.
Meromorphic braided category arising in quantum affine algebras, 
preprint math/9901003.
\bye